\documentclass[10pt]{amsart}
\usepackage{amssymb}
\newtheorem{thm}{Theorem}[section]
\newtheorem{cor}[thm]{Corollary}
\newtheorem{lem}[thm]{Lemma}
\newtheorem{prop}[thm]{Proposition}
\theoremstyle{definition}

\theoremstyle{remark}

\numberwithin{equation}{section}

\begin{document}
\title[Groups with certain elements]{finite groups with a certain number of elements
pairwise generating a non-nilpotent subgroup}%
\author{Alireza Abdollahi and Aliakbar Mohammadi Hassanabadi}%
\address{Department of Mathematics, University of
Isfahan, Isfahan 81744, and Institute for Studies in Theoretical
Physics and Mathematics, Iran} \email{a.abdollahi@sci.ui.ac.ir}
 \email{aamohaha@yahoo.com}
\thanks{This  research was in part supported by a grant from IPM}%
\thanks{Published in the {\sl Bulletin of the Iranian Mathematical
Sociey}, {\bf 30} No. 2 (2004), pp. 1-20.}
\begin{abstract}
Let $n>0$ be an integer and $\mathcal{X}$ be a class of groups.
 We say that a group $G$ satisfies the
condition $(\mathcal{X},n)$ whenever in every subset with $n+1$
elements of $G$ there exist distinct elements $x,y$ such that
$\left<x,y\right>$ is  in $\mathcal{X}$.  Let $\mathcal{N}$ and $
\mathcal{A}$ be the classes of nilpotent groups and abelian
groups, respectively. Here we prove that: (1)\;  If $G$ is a
finite semi-simple group satisfying the condition
$(\mathcal{N},n)$, then $|G|<c^{2[\log_{21}n]n^2} [\log_{21}n]!$,
for some constant $c$. (2)\; A finite insoluble group $G$
satisfies the condition $(\mathcal{N},21)$ if and only if
$\frac{G}{Z^*(G)}\cong A_5$, the alternating group of degree 5,
where $Z^*(G)$ is the hypercentre of $G$. (3)\;  A finite
non-nilpotent group $G$ satisfies the condition $(\mathcal{N}, 4)$
if and only if $\frac{G}{Z^*(G)}\cong S_3$, the symmetric group of
degree 3. (4)\; An insoluble group $G$ satisfies the condition
$(\mathcal{A},21)$ if and only if $G\cong Z(G)\times A_5$, where
$Z(G)$ is the centre of $G$. (5)\; If $d$ is the derived length of
a soluble group  satisfying the condition $(\mathcal{A},n)$,
 then $d=1$ if $n\in \{1,2\}$   and $d\leq 2n-3$ if
$n\geq 2$.\\
\end{abstract}
\maketitle
\section{Introduction and results}
Let $n>0$ be an integer and $\mathcal{X}$ be a class of groups.
 We say that a group $G$ satisfies the
condition $(\mathcal{X},n)$ whenever in every subset with $n+1$
elements of $G$ there exist distinct elements $x,y$ such that
$\left<x,y\right>$ is  in $\mathcal{X}$. If $\mathcal{X}$ is
subgroup-closed, then every
 group which  is the union of $n$ $\mathcal{X}$-subgroups satisfies the
condition $(\mathcal{X},n)$. Let $\mathcal{N}$ be the class of
nilpotent groups. Tomkinson in \cite{Tom} proved that if $G$ is a
finitely generated soluble group  satisfying the condition
$(\mathcal{N},n)$, then $|G/Z^*(G)|<n^{n^4}$, where $Z^*(G)$ is
the hypercentre of $G$.   This result gives  a bound for the size
of every
 finite soluble centerless group satisfying the condition $(\mathcal{N},n)$; on
the other hand, Endimioni in \cite{En} proved that if $n\leq 20$,
then every
 finite group satisfying the condition $(\mathcal{N},n)$ is soluble, and $A_5$, the
alternating group of degree $5$, satisfies the condition
$(\mathcal{N},21)$. Hence for $n\leq 20$ and all soluble groups,
we have a positive answer to the following question:\\ Does there
exist a bound (depending only on $n$) for  the
size of  every centerless finite group satisfying the condition $(\mathcal{N},n)$?\\
 Here we find a
bound for the size of finite semi-simple
 groups satisfying the condition $(\mathcal{N},n)$ and also for all
 finite centerless groups
satisfying the condition $(\mathcal{N},21)$. We also obtain a
characterization for $A_5$ (see Corollary \ref{c1.7}, below). The
main
results are \\

{\tt\sc Theorem  A.} {\sl Let $G$ be a finite semi-simple group
satisfying the condition $(\mathcal{N},n)$. Then
$|G|<c^{2[\log_{21}n]n^2}  [\log_{21}n]!$, for some constant $c$.}\\

{\tt\sc Theorem B.} {\sl Let $G$ be a finite insoluble group. Then
$G$ satisfies   the condition $(\mathcal{N},21)$  if and only if
$\frac{G}{Z^*(G)}\cong A_5$.\\}

In \cite{En} Endimioni proved that if $n\leq 3$, then every finite
group satisfying the condition $(\mathcal{N},n)$ is nilpotent, and
$S_3$, the symmetric group of degree 3, satisfies the condition
$(\mathcal{N},4)$. In fact, the only non-trivial finite centerless
group   satisfying the condition $(\mathcal{N},4)$ is $S_3$. In
section 2, we
investigate finite groups satisfying the condition $(\mathcal{N},4)$.\\

{\tt\sc Theorem C.} {\sl Let $G$ be a non-nilpotent finite group.
Then $G$  satisfies the condition $(\mathcal{N}, 4)$ if and
only if $\frac{G}{Z^*(G)} \cong S_3$. }\\

It follows from  Corollaries
  \ref{c1.8} and \ref{2.4} below that a finite
group satisfies the condition
 $(\mathcal{N},4)$  (respectively, $(\mathcal{N}, 21)$)
 if and only if it is the union of  4
 (respectively, 21)  nilpotent  subgroups. Another natural  question
is: ``For which positive integers $n$ is every finite group
satisfying the condition $(\mathcal{N},n)$ the union of $n$
nilpotent subgroups?" \\

In section 3, we investigate  (not necessarily finite) groups
satisfying the condition $(\mathcal{A},n)$, where $\mathcal{A}$ is
the class of abelian groups. Indeed, in a group  satisfying the
condition $(\mathcal{A},n)$, the largest  set of non-commuting
elements (or the largest set of elements in  which no two generate
an abelian  subgroup)
 has size at most $n$. By a result of B.H. Neumann  \cite{N}
a group satisfies the condition $(\mathcal{A},n)$ for some
$n\in\Bbb{N}$ if and only if it is  centre-by-finite. In fact,
Neumann  answered affirmatively the following question of P.
Erd\"os \cite{N}:
 Let $G$ be an infinite group. If there  is no
infinite subset of $G$ whose elements do not mutually commute, is
there then a finite bound on the cardinality of each such set of
elements? Neumann \cite{N} proved that a group  has the condition
of  Erd\"os's question if and only if it is centre-by-finite. This
result has  initiated a great deal of research towards the
determination of the structure of groups having some similar
properties (for example see
\cite{A1},\cite{A2},\cite{A3},\cite{AT1},\cite{AT2},\cite{D1},\cite{D2},
\cite{En2},\cite{Go},\cite{LW},\cite{LM},\cite{LMR},\cite{Tr}).\\
 Pyber in \cite{PY} gave a bound for the index of the
centre of a group   satisfying the condition $(\mathcal{A},n)$.
Here we  characterize insoluble groups satisfying the condition
$(\mathcal{A},21)$. Note that every group satisfying the condition
$(\mathcal{A},n)$ also satisfies the
condition $(\mathcal{N},n)$.\\

{\tt\sc Theorem D.} {\sl Let $G$ be an insoluble group. Then $G$
satisfies the condition $(\mathcal{A},21)$ if and only if $G\cong
Z(G)\times A_5$.\\}

We also obtain a result which is of independent interest,
namely,  the derived length of soluble groups satisfying the
condition
$(\mathcal{A},n)$ is bounded by a function  depending only on  $n$.\\

{\tt\sc Theorem E.} {\sl Let $G$ be a soluble group satisfying the
condition $(\mathcal{A},n)$ and let $d$ be the derived length of
$G$. Then $d=1$ if $n\in \{1,2\}$   and $d\leq 2n-3$ if
$n\geq 2$.}\\

\section{ Semi-simple groups satisfying the condition
$(\mathcal{N},n)$ and insoluble groups satisfying the condition
$(\mathcal{N},21)$}

Recall that a group $G$ is semi-simple if  $G$ has no non-trivial
normal abelian subgroups. If $G$ is a finite group then we call
the product of all minimal normal non-abelian subgroups of $G$ the
centerless CR-radical of $G$; it is a direct product of
non-abelian simple groups (see page 88 of \cite{R}).\\  We first
prove a result on the direct product of (not necessarily finite)
groups not satisfying the condition $(\mathcal{X},n)$, for a
certain class $\mathcal{X}$ of groups.
 This result may also be  useful in other investigations on
groups satisfying the condition $(\mathcal{X},n)$. For example, if
one can find a bound  depending only on $n$ for the size of finite
non-abelian simple groups satisfying the condition
$(\mathcal{X},n)$, then by the aid of  Lemma \ref{l1.1} below, it
is easy to see that there exists a bound depending only on $n$ for
the size of every semi-simple finite group satisfying the
condition $(\mathcal{X},n)$
(for instance  see Theorem A).\\
\begin{lem}\label{l1.1}
  Let $\mathcal{X}$ be a class of groups which is closed with respect
to homomorphic images. Suppose for $i\in\{1,\dots,t\}$ that $H_i$
is a group not satisfying the condition $(\mathcal{X},n_i)$. Then
$H_1\times\cdots\times H_t$ does not satisfy the condition
$(\mathcal{X},m)$, where $m=n_1+\cdots +n_t$.
 \end{lem}

\begin{proof} It suffices to show that if $H$ and $K$ are two
groups which do not satisfy $(\mathcal{X},n)$ and
$(\mathcal{X},m)$, respectively, then $H\times K$ does not satisfy
the condition $(\mathcal{X},n+m)$. By the hypothesis, there exist
$x_1,\dots,x_{n+1}$ in $H$ and $y_1,\dots,y_{m+1}$ in $K$ such
that $$\left<x_i,x_j\right>
 \not\in \mathcal{X} \;\;\text{for}\;\; 1\leq i<j\leq n+1
 \;\;\text{and}\;\;
 \left<y_k,y_l\right>\not\in\mathcal{X} \;\;\text{for}\;\; 1\leq
 k<l\leq m+1.$$
 Now it  is
easy to see that the subgroup generated by each pair of distinct
elements of the set
 $$\left\{(x_2,1),\dots,
(x_{n+1},1),(x_1,y_1),(x_1,y_2),\dots,(x_1,y_{m+1})\right\},$$
 does not have the property $\mathcal{X}$.
\end{proof}

Our next  lemma is about the direct product of finite groups not
satisfying $(\mathcal{N},n)$. For finite groups, this is  a better
result than Lemma \ref{l1.1}. \\
\begin{lem} \label{l1.2}   Suppose that $H_i$ is a finite
group not satisfying the condition $(\mathcal{N},n_i)$ for
$i\in\{1,\dots,t\}$. Then $H_1\times\cdots\times H_t$ does not
satisfy the condition $(\mathcal{N},m)$, where
$m=(n_1+1)\cdots(n_t+1)-1$.
\end{lem}

\begin{proof}
 By the hypothesis, for every $i\in\{1,\dots,t\}$ there exists a subset
$X_i$ in $H_i$ of size $n_i+1$  such that no pair of its distinct
elements generate a nilpotent subgroup.  Now we show that the
subgroup generated by each pair  of distinct elements of the set
$X=X_1\times\cdots\times X_t$ is not nilpotent. Let
$a=(a_1,\dots,a_t), b=(b_1,\dots,b_t)$ be two distinct  elements
of $X$. Then for some $i\in\{1,\dots, t\}$, $a_i\not
 =b_i$. Since $a_i, b_i \in X_i$, we have that $K:=\left<a_i,b_i\right>$ is not
 nilpotent. Since $K$ is a finite non-nilpotent group,
 it is not an Engel group by a result of Zorn (see Theorem
 12.3.4 of \cite{R}). Therefore there exist elements $x,y\in K$
 such that $[x,_ny]\not=1$ for all $n\in \Bbb{N}$. Suppose that
 $$x=a_i^{\delta_1} b_i^{\delta_2}\cdots
 a_i^{\delta_{r-1}}b_i^{\delta_r} \;\;\;\text{and}\;\;\; y=a_i^{\epsilon_1}
  b_i^{\epsilon_2}\cdots a_i^{\epsilon_{s-1}}b_i^{\epsilon_s}$$
  where $\delta_p,\epsilon_q\in\{0,1,-1\}$ for all $p\in\{1,\dots,r\}$
   and $q\in\{1,\dots, s\}$. Suppose, for a
  contradiction, that $\left<a,b\right>$ is nilpotent. Then there
  exists
   a positive  integer $m$ such that $[\bar{x},_m\bar{y}]=1$ where
   $$\bar{x}=a^{\delta_1} b^{\delta_2}\cdots
 a^{\delta_ {r-1}}b^{\delta_r} \;\;\;\text{and}\;\;\; \bar{y}=a^{\epsilon_1}
  b^{\epsilon_2}\cdots a^{\epsilon_{s-1}}b^{\epsilon_s}.$$
But $$[\bar{x},_m\bar{y}]=([x_1,_my_1],\dots, [x_t,_my_t])$$ where
$$x_j=a_j^{\delta_1} b_j^{\delta_2}\cdots
 a_j^{\delta_{r-1}}b_j^{\delta_r} \;\text{and}\; y_j=a_j^{\epsilon_1}
  b_j^{\epsilon_2}\cdots a_j^{\epsilon_{s-1}}b_j^{\epsilon_s}$$
for all $j\in\{1,\dots,t\}$.  Hence $[x,_my]=[x_i,_my_i]=1$, a
contradiction. This completes the proof.
\end{proof}
\begin{lem}\label{l1.3}  Let $M_1,\dots, M_m$ be non-abelian finite
simple groups. Then $M_1\times\cdots\times M_m$ does not satisfy
the condition $(\mathcal{N},21^m-1)$.
\end{lem}
\begin{proof} Since by  Proposition 2 of \cite{En}, $M_i$ does not
satisfy the condition $(\mathcal{N},20)$ for all
$i\in\{1,\dots,m\}$,
 the proof follows easily from  Lemma \ref{l1.2}.
 \end{proof}

Now we are ready to prove Theorem A.\\

{ \sc Proof of Theorem A.} Let $R$ be the centerless CR-radical of
$G$. Then $R$ is a direct product of a finite number $m$ of finite
non-abelian simple groups and $G$ is embedded in $Aut(R)$. Then by
Lemma \ref{l1.3}, we have  $21^m-1<n$ and so $m\leq [\log_{21}n]$.
On the other hand,
 since $Z(G)=1$, by Lemma 3.3 of \cite{Tom}
every prime divisor of $G$ is less than $n$. Thus by  Remark 5.5
of \cite{BGP}, there is a constant $c$
 such that   the order of every non-abelian simple section of $G$ is less
than $c^{n^2}$. Hence $|R|< c^{n^2[\log_{21}n]}$. Now using the
following well-known facts that: (a) \ for a finite simple group
$S$ we have $|Aut(S)|<|S|^2$ and (b) \ if $R$ is the product of
$m$ simple groups $S_i$, then $G$ acts on these  factors, the
quotient group is embeddable into $Sym(m)$ and the kernel $K$ of
the action is embeddable into the product of groups $Aut(S_i)$;
hence $|K|<|R|^2$. Thus $|G|<c^{2n^2[\log_{21}n]}
[\log_{21}n]!$.  \;\;\;$\Box$\\

Since in every finite group $G$, the quotient $G/Sol(G)$ is
semisimple, where $Sol(G)$ is the soluble radical (the largest
soluble normal subgroup) of $G$, we have

\begin{cor}\label{cor1}
Let $G$ be a finite group satisfying $(\mathcal{N},n)$. Then
$$|G/Sol(G)|<c^{2n^2[\log_{21}n]} [\log_{21}n]!$$ for some
constant $c$.
\end{cor}

Combining  the result of Tomkinson quoted in the introduction and
Corollary \ref{cor1}, we obtain as a further nice corollary that
in fact:
\begin{cor}
Let $G$ be a finite group satisfying $(\mathcal{N},n)$. Then
$$|G/F(G)|<n^{n^4}c^{2n^2[\log_{21}n]}[\log_{21}n]!$$ for some
constant $c$, where $F(G)$ is the largest nilpotent normal
subgroup of $G$.
\end{cor}
We need  the following proposition, which is of independent
interest, in the proof of Proposition \ref{p1.4}.
\begin{prop}\label{pppp}
Let $p$ be a prime number,  $n$ a positive integer  and $r$ and
$q$ be two odd prime numbers dividing respectively $p^n+1$ and
$p^n-1$. Then the number of Sylow $r$-subgroups (respectively,
$q$-subgroups) of $G=\text{PSL}(2,p^n)$ is $\frac{p^n(p^n-1)}{2}$
(respectively, $\frac{p^n(p^n+1)}{2}$). Also the intersection of
every two distinct Sylow $r$-subgroups or $q$-subgroups is
trivial.
\end{prop}
\begin{proof}
Our proof uses Theorems 8.3 and 8.4 in chapter II of \cite{H}.\\
Let $q$ be an odd prime dividing $p^n-1$ and let
$k=\text{gcd}(p^n-1,2)$. By Theorem 8.3 in Chapter II of \cite{H},
$\text{PSL}(2,p^n)$ possesses a cyclic subgroup $U$ of order
$u=\frac{p^n-1}{k}$ such that
\begin{enumerate}
\item The intersection of every two distinct conjugates of $U$ is
trivial.
\item For every non-trivial element $w$ of $U$, the normalizer
$N_G(\left<w\right>)$ of $\left<w\right>$ is a dihedral group of
order $2u$.
\end{enumerate}
Since $q$ is an odd prime number, $q$ divides $u$, and since
$|G|=\frac{p^n(p^n+1)(p^n-1)}{k}$, we have
$\text{gcd}(p^n(p^n+1),q)=1$. It follows that  any Sylow
$q$-subgroup of $U$ is also a Sylow $q$-subgroup of $G$ and each
of them is  cyclic. Therefore it follows from (2) that the number
of Sylow $q$-subgroups of $G$ is $\frac{p^n(p^n+1)}{2}$. Now (1)
implies that the intersection of every two distinct Sylow
$q$-subgroups of $G$ is trivial.\\
By a  similar argument the second statement of the proposition
 follows from the corresponding parts of  Theorem 8.4
in Chapter II of \cite{H}, namely that the group $G$  contains a
cyclic subgroup $K$ of order $s=\frac{p^n+1}{k}$ such that
\begin{enumerate}
\item The intersection of every two distinct conjugates of $K$ is
trivial.
\item For every non-trivial element $t$ of $K$, the normalizer
$N_G(\left<t\right>)$ of $\left<t\right>$ is a dihedral group of
order $2s$.
\end{enumerate}
\end{proof}
\begin{prop}\label{p1.4}
 The only non-abelian finite simple group satisfying the
condition $(\mathcal{N},21)$ is $A_5$.
\end{prop}
\begin{proof} Suppose, for a contradiction, that there exists
 a non-abelian finite  simple group satisfying the condition
 $(\mathcal{N},21)$
which is not isomorphic to $A_5$. Let $G$ be  such a group of
least order. Thus every proper non-abelian simple section of $G$
is isomorphic to $A_5$. Therefore by Proposition 3 of \cite{BR},
$G$ is isomorphic to one of the following:\\ \noindent
$\text{PSL}(2,2^p)$, $p=4$ or a prime;\\ $\text{PSL}(2,3^p)$,
$\text{PSL}(2,5^p)$, $p$ a prime;\\ $\text{PSL}(2,p)$, $p$ a prime
$\geq 7$;\\ $\text{PSL}(3,3)$, $\text{PSL}(3,5)$;\\
$\text{PSU}(3,4)$ (the projective special unitary group of degree
3 over the finite field of order $4^2$) or
\\ $Sz(2^p)$, $p$ an odd prime.\\ For each prime divisor $p$ of
$|G|$, let $\nu_p(G)$ be the number  of  all Sylow $p$-subgroups
of $G$. If $p$ is a prime number dividing $|G|$ such that the
intersection of any two distinct
Sylow $p$-Subgroups is trivial, then  by Lemma 3 of \cite{En}, $\nu_p(G)\leq 21$ (*).\\
 Now, for every prime
number $p$ and every integer  $n>0$, we have
$\nu_p(\text{PSL}(2,p^n))= p^n+1$ and the intersection of any two
distinct Sylow $p$-subgroups is trivial (see chapter II Theorem
8.2 (b),(c) of \cite{H}). Thus among the projective special
linear groups, we only need to investigate the following:
 $$\text{PSL}(2,3^2), \text{PSL}(2,8),
\text{PSL}(2,2^4), \text{PSL}(3,3), \text{PSL}(3,5),
\text{PSL}(2,p)$$ for $p\in\{7, 11, 13, 17, 19\}$. Now if in
 Proposition \ref{pppp}, we take $q=7$ for $\text{PSL}(2,8)$;
$q=5$ for $\text{PSL}(2,16)$; $r=5$ for $\text{PSL}(2,9)$; $q=3$
for $\text{PSL}(2,7)$, $\text{PSL}(2,13)$ and \text{PSL}(2,19);
and $r=3$ for $\text{PSL}(2,11)$ and $\text{PSL}(2,17)$;  we see,
by (*), that $G$ cannot be isomorphic with any of these groups.\\
Therefore we must consider  the groups $\text{PSL}(3,3),
\text{PSL}(3,5), \text{PSU}(3,4)$ or $Sz(2^p)$, $p$ an odd
prime.\\
$H:=\text{PSL}(3,3)$ has order $2^4\times 3^3\times 13$, so
$\nu_{13}(H)=1+13k$, for some $k>0$  and since $14$ does not divide $|H|$, $\nu_{13}(H)>26$.\\
$K:=\text{PSL}(3,5)$ has order $5^3\times 2^5 \times 3\times 31$,
so $\nu_{31}(K)=1+31k>21$ for some $k>0$. \\
$L:=\text{PSU}(3,4)$ has order $2^6\time 3\times 5^2\times 13$
(see Theorem 10.12(d) of chapter II in \cite{H} and note that $L$
is the projective special unitary group of degree 3 over the
finite field of order $4^2$). Therefore $\nu_{13}(L)=1+13k>21$ for
some $k>0$ and since $14$ does not divide $|L|$, $\nu_{13}(L)>26$. \\
$M:=\text{Sz}(2^p)$ ($p$ an odd prime) has order
$2^{2p}(2^p-1)(2^{2p}+1)$ and $\nu_2(M)= 2^{2p}+1\geq 65$ (see
Theorem 3.10 (and its proof) of chapter XI  in \cite{HB}). This
completes the proof by (*).
\end{proof}

\begin{lem}\label{l1.5}   $S_5$, the symmetric group of
degree 5, does not satisfy the condition $(\mathcal{N},21)$.
\end{lem}

\begin{proof} Every subgroup generated by a pair of  distinct
elements of 22-element subset\\
 \{(3,4,5), (2,3,4), (2,3,4,5),(1,4,5),
(2,3,5,4), (2,3,5), (2,4,5), (1,2,3),  (1,2,3,4),\\ (1,2,4,5,3),
(1,2,4,3,5), (1,2,5),(1,3,4), (1,3,4,5), (1,3,5), (1,3,2,4,5),
(1,4,2),\\(1,5,4,3,2), (1,5,3,2), (1,5,4,2),
 (1,5,2,4,3), (1,5,3,2,4)\}
  is not nilpotent. \end{proof}

{\sc Remark 1.} Here we state two properties of $A_5$ which we use
in the sequel.  Suppose that $P_1,\dots, P_{21}$ are all the Sylow
subgroups of $A_5$. Then \\ (i)\; For all $x_i\in
P_i\backslash\{1\}$ ($i=1,\dots,21$), the set
$\{x_1,\dots,x_{21}\}$ is a subset of $A_5$ such that no pair of
its distinct elements generate a nilpotent subgroup. (See the
proof of Proposition 2 of
\cite{En}).\\ (ii)\; $A_5=\cup_{i=1}^{21}P_i$.\\

We use the following fact in the sequel without any specific
reference.
 If $G$ is any group such that $G/Z_m(G)$ is nilpotent
for some integer $m\geq 0$, then  $G$ is nilpotent. For
$Z_n(\frac{G}{Z_m(G)})=\frac{G}{Z_m(G)}$ for some integer $n\geq
0$ and so by Theorem 5.1.11 (iv) of \cite{R}, we have
$Z_{m+n}(G)=G$, which implies that $G$ is nilpotent.\\

\begin{lem} \label{l1.6} Let $G$ be a finite insoluble group
satisfying the condition $(\mathcal{N},21)$   and let $S=Sol(G)$
be the soluble radical of $G$. Then $\frac{G}{S}\cong A_5$, and
for all $a\in S$ and for all $x\in G\backslash S$, the subgroup
$\left<a,x\right>$ is nilpotent. In particular, $Z^*(G)=Z^*(S)$.
\end{lem}
\begin{proof} Let $S$ be the soluble radical   of $G$ and
consider the semi-simple
 group $\overline{G}=G/S$. Let $\overline{R}$  be the
 centerless  CR-radical of $\overline{G}$. Then $\overline{R}$ is a direct
product of  non-abelian simple groups. Since $G$ is insoluble,
$\overline{R}$ is non-trivial. Now, by Lemma \ref{l1.3} and
Proposition \ref{p1.4}, $\overline{R}\cong A_5$. Since
$C_{\overline{G}}(\overline{R})=1$, we have $\overline{G}\cong
A_5$ or $S_5$. By Lemma \ref{l1.5}, $\overline{G}\cong A_5$.
 Now, let $Q_1,\dots,Q_{21}$ be the Sylow
subgroups of $G/S$.  For each $i\in\{1,\dots,21\}$, let $x_iS$ be
a non-trivial element of $Q_i$. Then, by Remark 1(i),
$\left<x_i,x_j\right>S\not\in \mathcal{N}$ and  so
$\left<x_i,x_j\right>\not\in \mathcal{N}$ for all distinct
$i,j\in\{1,\dots,21\}$. Now, fix $k\in\{1,\dots,21\}$ and for an
arbitrary element $a\in S$ consider the elements
$$x_k,x_1,\dots,x_{k-1},ax_k,x_{k+1},\dots,x_{21}.$$ For
$k,j\in\{1,\dots,21\}$ and $j\not=k$, $\left<ax_k,x_j\right>$ is
not nilpotent, since
$\left<ax_k,x_j\right>S=\left<x_k,x_j\right>S$. Since $G$
satisfies the condition $(\mathcal{N},21)$, the subgroup
$\left<x_k,ax_k\right>$ is nilpotent and hence so is
$\left<a,x_k\right>$ for all $k\in \{1,\dots,21\}$. On the other
hand, the union of the subgroups $Q_1,\dots,Q_{21}$ is $G/S$, by
Remark 1(ii), and so $\left<a,x\right>$ is nilpotent for all $x\in
G\backslash S$ and for all
$a\in S$.\\
 Since $S$ is finite, $Z^*(S)=Z_m(S)$ for some $m\in \Bbb{N}$. Now for all
$a\in Z_m(S)$ and for all $b\in S$, the subgroup
$T:=\left<a,b\right>$ is nilpotent, since $TZ_m(S)/Z_m(S)\cong
T/(T\cap Z_m(S))$ is cyclic and $T\cap Z_m(S)\leq Z_m(T)$. Thus
$\left<a,x\right>$ is nilpotent for all $a\in Z^*(S)$
 and for all $x\in G$. Since $G$ is finite, $a$ is a right Engel element for all
  $a\in
 Z^*(S)$ (see Theorem 12.3.7 of \cite{R})
  and so $Z^*(S)\leq Z^*(G)$. Hence $Z^*(S)=Z^*(G)$.
 This completes the proof.
 \end{proof}

{\sc Proof of Theorem B.} Suppose that $G$ satisfies the condition
$(\mathcal{N},21)$ and suppose, for a contradiction, that $G$ is a
counterexample of least order subject to $\frac{G}{Z^*(G)}\not
\cong A_5$. Let $S=Sol(G)$ be the soluble radical of $G$. We claim
that $Z(S)=1$. For if $Z(S)\not =1$ then $G/Z(S)$ is a finite
insoluble group satisfying the condition $(\mathcal{N},21)$ and
since $|\frac{G}{Z(S)}|<|G|$ and the soluble radical of $G/Z(S)$
is $S/Z(S)$, we have that the assertion of Theorem B is true for
the group $G/Z(S)$, i.e.
$$\frac{G/Z(S)}{Z^*(G/Z(S))}\cong A_5. \eqno{(*)}$$
Now Lemma \ref{l1.6}  implies that $Z^*(S/Z(S))=Z^*(G/Z(S))$. On
the other hand  $$Z^*(S/Z(S))=Z^*(S)/Z(S)=Z^*(G)/Z(S),$$ by Lemma
\ref{l1.6} (note that for a finite group $K$ we have
$Z^*(K)=Z_m(K)$ for some integer $m>0$). Thus it follows from
$(*)$ that
  $G/Z^*(G)\cong A_5$ which is a contradiction.  Hence
$Z(S)=1$, which implies that $Z^*(S)=1$.\\
 Now, let $x\in G\backslash S$ be such that $x^2\in S$. Thus for
all $b\in S$, we have $bx\in G\backslash S$ and $(bx)^2\in S$. By
Lemma \ref{l1.6}, $\left<bx,a\right>$ is nilpotent for all $a\in
S$, and so   also is $\left<(bx)^2,a\right>$. Therefore $(bx)^2$
is a right Engel element of $S$ and so $(bx)^2\in Z^*(S)=1$. Thus
for all $b\in S$, $(bx)^2=1$. Now, again by Lemma \ref{l1.6},
$\left<bx,x\right>=\left<b,x\right>$ is nilpotent and so is
$\left<b,x^2\right>$. Thus as before $x^2=1$. Therefore
$D:=\left<b,x\right>$ is a finite dihedral group which is
nilpotent and so $|D|$ is a power of 2 and $b$ is a 2-element.
Hence $S$ is a 2-group, and since $Z(S)=1$, we conclude that  $S$
must be trivial. Therefore, by Lemma \ref{l1.6}, $Z^*(G)=1$ and
$G/Z^*(G)=G/S\cong A_5$, a contradiction.  \\

Conversely, suppose that $\frac{G}{Z^*(G)}\cong A_5$. By Remark
1(ii), $$\frac{G}{Z^*(G)}=\bigcup_{i=1}^{21}\frac{P_i}{Z^*(G)},$$
where $\frac{P_1}{Z^*(G)},\dots, \frac{P_{21}}{Z^*(G)}$ are the
Sylow subgroups of $\frac{G}{Z^*(G)}$. But $G$ is finite, so
$Z^*(G)=Z_m(G)$ for some $m\in \Bbb{N}$. Since $Z_m(G)\leq
Z_m(P_i)$ for all $i\in\{1,\dots,21\}$ and $P_i/Z_m(G)$ is
nilpotent,  we conclude that each $P_i$ is nilpotent. Now the
proof is complete since $G=\cup_{i=1}^{21}P_i$. \;\;\;$\Box$\\

From Theorem B we have  a nice  characterization for $A_5$.\\

\begin{cor}\label{c1.7} The only finite  centerless
insoluble group satisfying the condition $(\mathcal{N},21)$ is
$A_5$.\end{cor}

Theorem B also gives us the following consequences.\\

\begin{cor}\label{c1.8}  A finite insoluble group
satisfies the condition $(\mathcal{N}, 21)$ if and only if it is
covered by 21 nilpotent subgroups.
\end{cor}

\begin{cor} \label{c1.9}  Let $G$ be a finite group satisfying
the condition $(\mathcal{N},21)$. If the centerless CR-radical of
$G$ is non-trivial, then $G\cong A_5 \times Z^*(G)$.\end{cor}

\begin{proof} Let $R$ be the centerless CR-radical of $G$. Then
$R$ is a non-trivial direct product of some non-abelian simple
groups and so by Lemma \ref{l1.3} and Proposition \ref{p1.4},
$R\cong A_5$. Since $R$ is simple, $R\cap Z^*(G) =1$. But, by
Theorem B, $|G|=|Z^*(G)||A_5|$, and so $G\cong A_5 \times Z^*(G)$.
\end{proof}

{\sc Remark 2.} We note that not every finite insoluble  group
satisfying the condition $(\mathcal{N} ,21)$ is  necessarily
isomorphic to a direct product as in Corollary \ref{c1.9}. For
example if $K:=SL(2,5)$ then $\frac{K}{Z(K)}\cong A_5$ and so $K$
satisfies the condition $(\mathcal{N},21)$, by Theorem B. However
we conjecture that every finite insoluble group satisfying the
condition $(\mathcal{N},21)$ is a direct product of a nilpotent
group and a group   isomorphic to either $A_5$ or $SL(2,5)$.

\section{ Finite groups satisfying the condition $(\mathcal{N},4)$}

In this section, we investigate finite groups satisfying the
condition $(\mathcal{N},4)$,  and give the proof of Theorem C.\\

\begin{lem} \label{2.1} Let $G$ be a finite $\{2,3\}$-group. If $G$
satisfies the condition $(\mathcal{N},4)$, then $G$ is
2-nilpotent.\end{lem}

\begin{proof} Suppose that $G$ is a counterexample of  least
order. Thus by a result of It$\hat{\text{o}}$ (see Theorem 5.4 on
page 434 of \cite{H}), $G$ is a minimal non-nilpotent group  and
$G$ has  a unique Sylow $2$-subgroup $P$ and a cyclic Sylow
$3$-subgroup $Q$ such that $\Phi(Q)\leq Z(G)$ and   $\Phi(P)\leq
Z(G)$ (see Theorem 5.2 on page 281 of \cite{H}). If $Z(G)\not =1$
then $G/Z(G)$ is nilpotent and so $G$ is nilpotent, a
contradiction. Thus $Z(G)=1$
 and so $|Q|=3$  and $P$ is an elementary abelian  2-group. Let
$Q=\left<a\right>$. Then $C_P(a)\leq Z(G)$,  and so $C_P(a)=1$. On
the other hand by Lemma 3.4 of \cite{Tom}, $|P:C_P(a)|\leq 4$ and
so $|P|\leq 4$. If $|P|=4$ then $G\cong A_4$, the alternating
group of degree 4. But $A_4$ does not satisfy the condition
$(\mathcal{N},4)$; thus $|P|=2$. Therefore $G\cong S_3$, a
contradiction, since $S_3$ is 2-nilpotent. This completes the
proof. \end{proof}

\begin{lem} \label{2.2}  Let $G=RX$ be an extension of an elementary
abelian 3-group $R$ by an abelian 2-group $X$ such that $X$ acts
faithfully on $R$ and $R=[R,X]$.  If $G$ satisfies the condition
$(\mathcal{N},4)$, then $|X|\leq 2$ and $|R|\leq 3$.
\end{lem}

\begin{proof} The proof  follows from the argument of  Lemma 3.7 of
\cite{Tom}.
\end{proof}

We are now ready to give a proof for  Theorem C,  the outline  of
which is in fact a refinement of that of  Theorem C in \cite{Tom}
for $n=4$.\\

{\sc Proof of  Theorem C.} Suppose that $G$  satisfies the
condition $(\mathcal{N},4)$. By factoring out $Z^*(G)$, we may
assume that $G$ is a finite non-trivial group with trivial centre
satisfying the condition $(\mathcal{N}, 4)$.
 We note that   $G$ is a
$\{2,3\}$-group by Lemma 3.3 of \cite{Tom}. \\ Let $H_p/O_{p'}(G)$
be the hypercentre of $G/O_{p'}(G)$, for $p=2,3$. Then, since $G$
is finite, there is a positive  integer $m$ such that $[H_p,
_mG]\leq O_{p'}(G)$ for $p=2,3$. Hence $$[H_2 \cap H_3, _mG]\leq
O_{2'}(G) \cap O_{3'}(G)=1$$ and so $H_2 \cap H_3 \leq Z^*(G)=1$.
But $O_{2'}(G)=O_3(G)$ and by Lemma \ref{2.1}, is the unique Sylow
3-subgroup of $G$. Thus $G/O_{2'}(G)$ is a 2-group and so $G=H_2$.
Therefore $H_3=1$ and so $O_{2}(G)=1$.  Hence
$P=\text{Fitt}(G)=O_3(G)$. Let $\overline{G}=G/\Phi(P)$ and
$\overline{P}=P/\Phi(P)$, thus $\overline{G}/\overline{P}$ acts
faithfully on the $GF(3)$-vector space $\overline{P}$ (see
\cite{Gr}, Theorem 6.3.4). We note that $\overline{P}$  is an
elementary abelian normal 3-subgroup of $\overline{G}$,  that
$\overline{P}=O_3(\overline{G})$,  and that
$C_{\overline{G}}(\overline{P})=\overline{P}$. Let
$Q/\overline{P}$ be the socle of $\overline{G}/\overline{P}$, so
that $Q/\overline{P}$ is an abelian 2-subgroup. We may write
$Q=\overline{P}X$, where $X$ is an abelian 2-subgroup of $Q$. Let
$R=[\overline{P},Q]$, so that $\overline{P}=R\times
C_{\overline{P}}(Q)$. If $C=C_{\overline{G}}(R)$ then $C\cap Q$
centralizes $R \times C_{\overline{P}}(Q)=\overline{P}$ and so
$C\cap Q =\overline{P}$. It follows that
$C_{\overline{G}}(R)=\overline{P}$ and so
$\overline{G}/\overline{P}$ acts faithfully on $R$. Now $R$ and
$X$ satisfy the conditions of Lemma \ref{2.2} and so  $|R|\leq 3$.
Since $\overline{G}/\overline{P}$ acts faithfully on $R$,
 the order of $G/P$ is no more than 2.
  Let $T$ be a Sylow 2-subgroup of $G$; then $|T|\leq 2$ and hence $T$
is cyclic and by Lemma 3.4 of \cite{Tom}, $|P:C_P(T)|\leq 3$. Now,
we have $[C_P(T), _mG]=[C_P(T),_mP]=1$ for some $m\in\Bbb{N}$.
Thus $C_P(T)\leq Z^*(G)=1$ and so $|G|=|T||P|\leq 2\times 3=6$.
Therefore $G\cong S_3$. \\
Conversely, suppose that $G/Z^*(G)\cong S_3$. Since $S_3$ is
covered by 4 abelian subgroups, $G$ is also covered by 4 nilpotent
subgroups. This completes the proof. \;\;\;$\Box$\\

\begin{cor}\label{2.3}  Every finite group satisfying the
condition $(\mathcal{N},4)$ is supersoluble. The alternating
group $A_4$ satisfies the condition $(\mathcal{N},5)$.
\end{cor}

\begin{proof} Let $G$ be a finite group satisfying the condition
$(\mathcal{N},4)$. By Proposition 1 of \cite{En}, $G= H\times K$,
where $H$ is a nilpotent $\{2,3\}'$-group and $K$ is a
$\{2,3\}$-group. If $K$ is nilpotent, then there is nothing to
prove. Assume that $K$ is not nilpotent. By Theorem C,
$K/Z^*(K)\cong S_3$ and so $K$ is supersoluble. Thus $G$ is also a
supersoluble group.\\ The group $A_4$ is the union of its five
Sylow subgroups , so $A_4$ satisfies the condition
$(\mathcal{N},5)$.
\end{proof}

\begin{cor}\label{2.4}  A finite group satisfies the
condition $(\mathcal{N},4)$ if and only if it is the union of four
nilpotent  subgroups.
\end{cor}

\begin{proof} Let $G$ be a finite group satisfying the condition
$(\mathcal{N},4)$. Then by  Theorem C,  $G/Z^*(G)$ is the union of
4 nilpotent  subgroups and hence so is $G$. The converse is clear.
\end{proof}

\section{ Finite groups satisfying the condition $(\mathcal{A},n)$}

 Now suppose that $\mathcal{A}$ is
the class of abelian groups. Then every group satisfying the
condition $(\mathcal{A},n)$ also satisfies the condition
$(\mathcal{N},n)$. The converse is not true, since, as we have
observed already,  $SL(2,5)$ satisfies the condition
$(\mathcal{N},21)$. However $SL(2,5)$ does not satisfy
the condition $(\mathcal{A},21)$.\\

\begin{lem} \label{3.1}  $SL(2,5)$ does not satisfy the
condition $(\mathcal{A},21)$.
\end{lem}

\begin{proof} Let $P_1, \dots, P_5$ be the Sylow 2-subgroups of $SL(2,5)$,
$Q_1,\dots, Q_{10}$ the Sylow 3-subgroups of $SL(2,5)$, and
$R_1,\dots, R_6$ the Sylow 5-subgroups of $SL(2,5)$. We note for
each $i=1,\dots,5$ that $P_i$ is a quaternion group of order 8 and
$Z(P_i)=Z(SL(2,5))$
 (see, for example, Theorem 8.10 in chapter II
of \cite{H}). Let $x_i\in P_i\backslash Z(P_i)$ ($i=1,\dots,5$),
$y_{j}\in Q_j\backslash\{1\}$ ($j=1,\dots,10$) and $z_{k}\in
R_{k}\backslash\{1\}$ ($k=1,\dots,6$). Then since
 $\frac{SL(2,5)}{Z(SL(2,5))}\cong A_5$,  it follows from Remark  1(i) following Lemma \ref{l1.5} that
no two distinct elements of the set $$\{x_1,\dots,
x_5,y_1,\dots,y_{10},z_1,\dots, z_6\}$$ commute.
  Now since $P_1$ is a quaternion group of order 8   and
$x_1\in P_1\backslash Z(P_1)$, there exists an element $x\in
P_1\backslash Z(P_1)$ such that $x_1x\not=xx_1$. On the other
hand, as above,  no two distinct elements in $$\{x,x_2\dots,
x_5,y_1,\dots,y_{10},z_1,\dots, z_6\}$$  commute. Therefore no two
distinct elements in the set $$\{x,x_1,\dots,
x_5,y_1,\dots,y_{10},z_1,\dots, z_6\}$$  commute, which completes
the proof.
\end{proof}
\begin{lem} \label{3.2}  Let $G$ be a finite group satisfying the
condition $(\mathcal{A},21)$. If there exists a central subgroup
$B$ of $G$ of order no more than 2 such that $G/B\cong A_5$, then
$G\cong B\times A_5$.
\end{lem}
\begin{proof}   Since $G/B\cong A_5$ it follows that
$G=G'B$ and $G'/(B\cap G')\cong A_5$. Therefore if $G'\cap B=1$
then the proof is complete. So suppose, for a contradiction, that
$G'\cap B\not=1$. Thus $|B|=2$. According to the Universal
Coefficients Theorem (see Theorem 11.4.18 of \cite{R}) the
central extension $B\rightarrowtail G\twoheadrightarrow G/B$
determines a homomorphism $\delta: M(\frac{G}{B})\rightarrow B$ so
that $\text{Im} \delta =G' \cap B$,  where  $M(\frac{G}{B})$ is
the Shur multiplicator of $\frac{G}{B}$ (see for example
Exercise  10 on page 354 of \cite{R}). But we know that the Shur
multiplicator of the alternating group $A_5$ is $\mathbb{Z}_2$.
Hence $G'\cap B=B$ and so $B\leq G'$. It follows that $G$ is a
perfect group of order 120. But it is well-known that the only
perfect group of order $120$ is $\text{SL}(2,5)$. Now Lemma
\ref{3.1} gives a contradiction and the proof is complete.
\end{proof}
We need the following lemma in the proof of Theorem D.

\begin{lem} \label{3.3}  Let $G$ be a group satisfying the condition
$(\mathcal{A},n)$ ($n>1$). Then for any  normal non-abelian
subgroup $N$ of $G$,  the quotient $G/N$ satisfies the condition
$(\mathcal{A},n-1)$.
\end{lem}

\begin{proof} Suppose, for a contradiction, that $G/N\not\in
(\mathcal{A},n-1)$. Then there exist elements $x_1,\dots,x_n$ in
$G$ such that $[x_i,x_j]\not\in N$ for all distinct
$i,j\in\{1,\dots,n\}$ $(*)$. Let $a,b$ be two distinct arbitrary
elements of $N$ and consider the subset
$X=\{ax_1,\dots,ax_n,bx_1\}$. By the hypothesis, there exist two
distinct commuting elements in $X$. But, by $(*)$, the only
commuting pair of elements of $X$ are $bx_1$ and $ax_1$. Therefore
for all $a,b\in N$, we have $ax_1b=bx_1a$ $(**)$ and in particular
for $b=1$,  we have $ax_1=x_1a$ for all $a\in N$. Thus for all $x,
y\in N$ we have
$$xyx_1=xx_1y=yx_1x=yxx_1$$ (the middle equality follows from
$(**)$) and so $xy=yx$. Hence $N$ is abelian, a contradiction.
\end{proof}

{\sc Proof of Theorem D.} Suppose that  $G\cong Z(G) \times A_5$.
Then $G$ is covered by 21  abelian subgroups as $A_5$ has this
property, by Remark 1(ii) following Lemma \ref{l1.5}.\\

Now, suppose that $G$ satisfies the condition $(\mathcal{A},21)$.
Then by a famous Theorem of B. H. Neumann \cite{N}, $G/Z(G)$ is
finite. Thus, by Theorem B,  $$\frac{G/Z(G)}{Z^*(G/Z(G))}\cong
G/Z^*(G)\cong A_5.$$ If $H:=Z^*(G)$ is not abelian, then Lemma
\ref{3.3} shows that $A_5$ satisfies the condition
$(\mathcal{A},20)$, which contradicts  Proposition 2 of \cite{En}.
Thus $H$ is abelian;   we show that in fact  $H=Z(G)$. To prove
this let $Q_1,\dots,Q_{21}$ be the Sylow subgroups of
$\overline{G}:=G/H$. For each $i\in\{1,\dots,21\}$, let $x_iH$ be
a non-trivial element of $Q_i$. Then $[x_i,x_j]\not\in H$ and so
$[x_i,x_j]\not=1$ for all distinct $i,j\in\{1,\dots,21\}$, by
Remark 1(i) following Lemma \ref{l1.5}. Now, fix
$k\in\{1,\dots,21\}$ and consider the elements
$$x_k,x_1,\dots,x_{k-1},ax_k,x_{k+1},\dots,x_{21},$$ for an
arbitrary element $a\in H$. Then for $j\in\{1,\dots,21\}$ and
$j\not=k$, we have $[ax_k,x_j]\not=1$, since
$[ax_k,x_j]H=[x_k,x_j]H$. Since $G$ satisfies the condition
$(\mathcal{A},21)$, $[x_k,ax_k]=1$ and so $[a,x_k]=1$ for all
$k\in \{1,\dots,21\}$. Since the union of $Q_1,\dots,Q_{21}$ is
$\overline{G}$, by Remark 1(ii) following Lemma \ref{l1.5},  we
have $[a,x]=1$ for all
$x\in G\backslash H$ and for all $a\in H$. Therefore $H=Z(G)$.\\
Now by the same argument as in Lemma \ref{3.2}, considering the
central extension $Z(G)=H\rightarrowtail G\twoheadrightarrow
\overline{G}$, we have that $K=G' \cap Z(G)$ is of order no more
than $2$, $G=G'Z(G)$ and $G'/K\cong A_5$. Thus Lemma \ref{3.2}
implies that there is a subgroup $L$ of $G'$ such that $G'=K\times
L$ and $L\cong A_5$. Therefore $G=G'Z(G)=LKZ(G)=LZ(G)$ and it is
clear that $L \cap Z(G)=1$. Therefore $G=L\times Z(G)\cong A_5
\times Z(G)$.
 \;\;\;$\Box$\\

We end this paper by proving Theorem E.\\

{\sc Proof of Theorem E.} We first prove that if $n=2$, then $G$
is abelian. Consider two distinct elements $x,y\in G$. Then
$X=\{x,y,xy\}$ is a subset of size 3. Thus by the hypothesis two
distinct elements of X commute. But commutativity of each pair of
distinct elements of $X$ implies the commutativity of $x$ and $y$.
Hence $G$ is abelian.\\ Now suppose that $n\geq 2$ and use
induction on $n$. If $n=2$ then $G$ is abelian and $d=1$. So let
$n>2$. Then   $2< 2n-3$. Thus we may assume that $d>2$. Therefore
$G^{d-2}$ is not abelian and so $G/G^{d-2}$ satisfies the
condition $(\mathcal{A},n-1)$ by Lemma
 \ref{3.3}. Thus by
induction the derived length of $G/G^{d-2}$ is at most $2(n-1)-3$
and so
$d-2\leq 2(n-1)-3$. Hence $d\leq 2n-3$. \;\;\;$\Box$\\

{\bf Acknowledgements.} The authors would like to thank the
referees for their careful considerations and valuable
suggestions.

\end{document}